\documentclass[12pt]{asl}

\usepackage{amssymb}
\usepackage{enumerate}
\usepackage{mathrsfs}

\title{Borel reducibility and classification of von Neumann algebras}

\author{Rom\'an Sasyk}
\revauthor{Sasyk, Rom\'an}

\address{
Instituto de Ciencias\\
Universidad Nacional de General Sarmiento\\
J. M. Gutierrez 1150\\
(1613) Los Polvorines, Argentina}
\address{
Instituto Argentino  de Matem\'aticas-CONICET\\
Saavedra 15, Piso 3
(1083), Buenos Aires, Argentina}

\email{rsasyk@ungs.edu.ar}

\author{Asger T\" ornquist}
\revauthor{T\" ornquist, Asger}

\address{Kurt G\"odel Research Center\\
University of Vienna\\
W\"ahringer Strasse 25\\
1090 Vienna, Austria}

\email{asger@logic.univie.ac.at}

\thanks{Work on this paper was initiated when R. Sasyk was a post-doc at the University
of Ottawa and A. T\"ornquist was a post-doc at the University of
Toronto, and we thank both institutions for their support. We also
wish to extend our gratitude towards Ilijas Farah, whose many
suggestions and comments have vastly improved this paper. A. T\"
ornquist was supported in part by Danish Natural Science Research
Council grant no. 272-06-0211.}

     \DeclareMathOperator{\Aut}{Aut}

    \DeclareMathOperator{\Mod}{Mod}

    \DeclareMathOperator{\im}{im}
    \DeclareMathOperator{\Char}{Char}
    \DeclareMathOperator{\Ext}{Ext}
    \DeclareMathOperator{\vN}{vN}
    \DeclareMathOperator{\SL}{SL}
    \DeclareMathOperator{\HT}{HT}
    \DeclareMathOperator{\I}{I}
    \DeclareMathOperator{\II}{II}
    \DeclareMathOperator{\III}{III}

    \DeclareMathOperator{\gp}{\bf GP}

    \DeclareMathOperator{\itpf1}{ITPFI}
    
    \DeclareMathOperator{\Tr}{Tr}
    \DeclareMathOperator{\pt}{T}
    \newcommand{\actson}{\curvearrowright}

\def\R{{\mathbb R}}
\def\C{{\mathbb C}}

\def\N{{\mathbb N}}
\def\Z{{\mathbb Z}}

\def\F{{\mathbb F}}

\def\T{{\mathbb T}}

\newcounter{mypar}
\newcounter{thmcounter}

\begin{document}

\maketitle

\begin{abstract}
We announce some new results regarding the classification problem
for separable von Neumann algebras. Our results are obtained by applying the notion
of Borel reducibility and Hjorth's theory of turbulence to the
isomorphism relation for separable von Neumann algebras.
\end{abstract}

\section{Introduction}

Let $\mathcal H$ be an infinite dimensional separable complex
Hilbert space and denote by $\mathcal B(\mathcal H)$ the space of
bounded operators on $\mathcal H$, which we give the weak topology.
A separable von Neumann algebra is a weakly closed self-adjoint
subalgebra of $\mathcal B(\mathcal H)$. A von Neumann algebra is
called a {\it factor} if its center consists of the scalar multiples
of the identity. The factors make up the building blocks of von
Neumann algebra theory: Any von Neumann algebra can be represented
as a direct integral of factors (see \cite[III.1.6]{blackadar}). A
central problem in the theory of von Neumann algebras is to classify
factors up to isomorphism (see \cite{connes}). The first steps towards a classification
were obtained by Murray and von Neumann, \cite{murvn}, when they
introduced the notion of the {\it type} of a von Neumann algebra and
gave examples of factors in each of the classes. Another major advance towards classifying factors
was A. Connes' thesis \cite{connes2} where he further extended the notion of type to split the type $\III$ case in the subtypes
 $\III_\lambda$, $0\leq\lambda\leq 1$.

Denote by $\vN(\mathcal H)$ the set of von Neumann algebras on
$\mathcal H$. E. Effros \cite{effros1}, \cite{effros2} defined a
Borel structure on $\vN(\mathcal H)$ and proved that this structure
is standard and that the set of factors $\mathscr F(\mathcal H)$ is
Borel. One of Effros's goals  was to show that the set of
isomorphism classes of factors endowed with the quotient Borel
structure was not standard, which would imply the existence of
uncountably many non-isomorphic von Neumann factors. At that time
only few examples were known of non-isomorphic infinite dimensional
factors, thus a solution of this problem was of major importance.
The existence of a continuum of non-isomorphic factors was finally
shown by Powers in 1967 for the type $\III$ case and by McDuff in
1969 for the type $\II_1$ case. In 1971 J. Woods solved Effros's
problem: The isomorphism relation for factors is not smooth. In
modern terminology, what Woods proved was that $E_0$, the
equivalence relation on $2^\N$ of eventual equality, is Borel
reducible to isomorphism of so-called ITPFI factors (see \S 4.2 or
\cite[III.3.1]{blackadar} for a definition).

Until recently, Woods's result was one of the few theorems of its
kind in the study of von Neumann algebras. For example, it remained
an open problem to show that isomorphism of factors of type $\II_1$
is not smooth. In a forthcoming paper \cite{asro} we apply the
notion of Borel reducibility from descriptive set theory to obtain
information about the classification problem for separable von
Neumann factors.

Let us briefly recall the key notions surrounding Borel
reducibility, but otherwise refer to the excellent introduction in
\cite{hjorth}, or the survey \cite{kechris2}. If $E$ and $F$ are
equivalence relations on standard Borel spaces $X$ and $Y$,
respectively, then we say that $E$ is {\it Borel reducible} to $F$,
written $E\leq_B F$, if there is a Borel $f:X\to Y$ such that
$$
xEy\iff f(x) F f(y).
$$
Thus if $E\leq_B F$ then the points of $X$ can be classified up to
$E$ equivalence by a Borel assignment of invariants that we may
think of as $F$-equivalence classes. $E$ is {\it smooth} if it is
Borel reducible to the equality relation on $\R$. While smoothness
is desirable, it is most often too much to ask for: As mentioned
above, $E_0$ is not smooth. A more generous class of invariants
which seem natural to consider are countable groups, graphs, fields,
or other countable structures, considered up to isomorphism. Thus,
following \cite{hjorth}, we will say that an equivalence relation
$E$ is {\it classifiable by countable structures} if there is a
countable language $\mathcal L$ such that
$E\leq_B\simeq^{\Mod(\mathcal L)}$, where $\simeq^{\Mod(\mathcal
L)}$ denotes isomorphism in $\Mod(\mathcal L)$, the Polish space of
countable models of $\mathcal L$ with universe $\N$. We note that
$E_0$ may be seen to be classifiable by countable structures.

It turns out that even allowing this more generous invariant is not
enough: There are many natural classification problems in
mathematics where countable structures do not suffice as invariants.
Hjorth conceived of his theory of {\it turbulence} (see again
\cite{hjorth}) as a general tool to prove that various equivalence
relations are not classifiable by countable structures. One of the
early applications of this theory was due to Foreman and Weiss \cite{foremanweiss}, who
showed that the measure preserving ergodic transformations on the
unit interval are not classifiable, up to conjugacy, by countable
structures. Subsequently, similar results have been achieved for the
weaker notion of orbit equivalence (see \S 2) of measure preserving
actions of non-amenable groups, see \cite{torn2}, \cite{ikt}.

Our main results are

\newtheorem{theorem}{Theorem}
\begin{theorem}
The isomorphism relation for factors of type $\II_1$,
$\II_\infty$ and $\III_\lambda$, $0\leq \lambda\leq 1$ is not
classifiable by countable structures.
\end{theorem}

\begin{theorem}
If $\mathcal L$ is a countable language, then $\simeq^{\Mod(\mathcal
L)}<_B\simeq^{\mathscr F_{\II_1}}$, where $\simeq^{\mathscr
F_{\II_1}}$ denotes the isomorphism relation for factors of type
$\II_1$.
\end{theorem}

Since it is known that the isomorphism relation for countable
graphs, say, is complete analytic, we obtain the following as a
consequence of Theorem 2:

\begin{theorem}
The isomorphism relation $\simeq^{\mathscr F_{\II_1}}$ of factors of
type $\II_1$ is a complete analytic subset of $\mathscr
F_{\II_1}\times\mathscr F_{\II_1}$, where $\mathscr F_{\II_1}$
denotes the Borel set of $\II_1$ factors.
\end{theorem}

The proofs of these results rely heavily on results obtained by
Popa's novel ``deformation rigidity techniques'', in particular on
the class of $\HT$ factors (discussed below) introduced in
\cite{popa}, as well as Hjorth's theory of turbulence. In this paper
we will first in \S 2 give some background regarding von Neumann
algebra factors, in particular regarding the {\it group-measure
space construction}, which plays the starring role in all the proofs
above. In \S 3 we give an outline of the proofs. In \S 4 we discuss
some of the open problems and questions that remain.

\section{Von Neumann algebras}

A separable von Neumann Algebra is a weakly closed self-adjoint
algebra of operators on a separable complex Hilbert space. A von
Neumann algebra is called a factor if its center only consists of
the scalar multiples of the identity. A von Neumann algebra $N$ is
said to be {\it finite} if it admits a finite faithful normal
tracial state, i.e. a linear functional $\tau:N\to \C$ such that:
$\tau(x^*x)\geq 0$, $\tau(x^*x)=0$ iff $x=0$, $\tau(1)=1$\,,
$\tau(xy)=\tau(yx)$ and the unit ball of $N$ is complete with
respect to the norm given by the trace $\|x\|_\tau=\tau(x^*x)$. If
such a trace exists it need not to be unique, however, a fundamental
fact is that if a finite von Neumann algebra is also a factor, then
it has a unique trace.

Some basic examples of finite von Neumann algebras are:
\begin{enumerate}
\item $L^\infty(X,\mu)$, the set of essentially bounded measurable functions on a standard
Borel probability space $(X,\mu)$. They act by multiplication on $L^2(X,\mu)$.
Here the trace is given by the integral. The Borel functional calculus yields that
any separable Abelian finite von Neumann algebra is of this form. \item
$M_n(\C)$, the set of $n\times n$ complex matrices with the
normalized trace $\Tr_{M_{n}(\C)}$. Any finite dimensional von
Neumann factor is of this form.
\item $\bigoplus_{i=1}^k M_{n_i}(\C)$ with the trace given
by $\sum_{i=1}^kc_i\Tr_{M_{n_i}(\C)}$; $c_i>0$, $\sum_{i=1}^k
c_i=1$. Moreover any finite dimensional von Neumann algebra is of
this form.
\end{enumerate}

Von Neumann algebras are categorized into {\it types} according to
the behavior of the lattice of projections. The types are called
$\I$, $\II_1,\II_\infty$ and $\III_\lambda, 0\leq\lambda\leq 1$ ,
see \cite[III.1.4]{blackadar} or \cite{ringrose} for an
introduction. A von Neumann algebra is of type $\II_1$ if it is
finite and it doesn't have minimal projections (a projection $p$ in
a von Neumann algebra $M\subseteq\mathcal B(\mathcal H)$ is said to
minimal if there is no projection $q\in M$ with $0<q<p$, where $q<p$
means $\im(q)\varsubsetneq\im(p)$.) If $M$ is a $\II_1$ factor and
$\tau$ is the normalized trace on $M$ then
$$
\{\tau(p): p\in M \text{ is a projection}\}=[0,1].
$$

Factors of type $\II_\infty$ are of the form
$$
M\otimes \mathcal B(\mathcal H)
$$
where $M$ is a factor of type $\II_1$ and $\mathcal H$ is an
infinite dimensional Hilbert space. A factor of type $\II_\infty$
has a semifinite trace which is unique up to scaling.

Murray and von Neumann already exhibited examples of factors of type
$\II_1$. This was done using two fundamental constructions, the {\it
group von Neumann algebra} and the {\it group-measure space
construction}. Since both constructions still play a preponderant
role in the theory and they are at the core of some of our
arguments, we give here a more or less detailed account of how they
are constructed.

\subsection{The group von Neumann algebra.} Let $G$ be an infinite
discrete group.  $\ell^2(G)$ is the infinite dimensional Hilbert
space with orthonormal basis $\{\xi_g:g\in G\}$. The group $G$ acts
on $\ell^2(G)$ by the left regular representation
$u_g(\xi_h)=\xi_{gh}$. The {\it group von Neumann algebra} $L(G)$ is
the von Neumann algebra generated  by the unitary operators $u_g\in
\mathcal B(\ell^2(G))$, that is, $L(G)=\overline{\langle u_g: g\in
G\rangle}^{wo}$, the closure in the weak operator topology of the
algebra generated by the $u_g$. The trace is given by
$\tau(u_g)=\langle u_g(\xi_e),\xi_e\rangle$, where $e$ denotes the
identity of $G$. It is easy to show that $L(G)$ is a factor iff the
group $G$ is ICC (infinite conjugacy classes) i.e. for each $g\in
G\backslash \{e\}$ the conjugacy class $\{hgh^{-1}: h\in G\}$ is
infinite.

Let $\gp$ denote the Polish space of all countable groups with universe $\N$, and consider the equivalence relation $\sim_{\vN}$ in $\gp$ defined by
$$
G\sim_{\vN} H\iff L(G) \text{ is isomorphic to } L(H).
$$
We do not know how complex this equivalence relation is (cf. Problem 7 below.) Two outstanding (yet seemingly unrelated) open problems in the
theory are concerned with this equivalence relation:

(a) Is it true that $\F_n\not\sim_{\vN} \F_m$ when $n\neq m$, $n,m\geq 2$? That is, when is $L(\F_n)$  isomorphic to $L(\F_m)$? Here $\F_n$ denotes the free group on $n$ generators.
Free probability was first envisioned by Voiculescu as an attempt to
solve this problem (see \cite{voiculescu},
\cite{VoiculescuDykemaNica} for an overview.)

(b) If $G$ and $H$ are countably infinite ICC property (T) groups, does $G\sim_{\vN} H$ imply that $G$ is isomorphic to $H$? This problem is known as Connes's conjecture\footnote{See
\cite{behava} for more information regarding Kazhdan's property
(T).}.

The first appearance of property
(T) in the context of operator algebras is Connes result
\cite{connesT} stating that the fundamental group of the group von
Neumann algebra of an ICC property (T) group is countable. Recall
that the fundamental group of a $\II_1$  factor $M$ is defined as
$$
F(M) = \{\tau(p)/\tau(q) | pMp \simeq qMq\}
$$
where $p,q\in M$ are non-zero projections and $\tau$ denotes the
trace on $M$. The fundamental group is a subgroup of $\R_{>0}$. As a
consequence of his work on $\HT$ factors, Popa gave in \cite{popa}
the first example of a type $\II_1$ factor with trivial fundamental
group, solving a longstanding problem in the theory. Going back to
Connes's conjecture, in \cite{popa2} Popa gave what may be seen as a
partial affirmative answer to the conjecture but for actions of
property (T) groups. In \S 3 we explain some aspects of Popa's work
that are pertinent to our results while for a more thorough
introduction to Popa's theory and its many applications we refer the
reader to the survey papers \cite{connesBourbaki},  \cite{popaICM}
and \cite{vaesBourbaki}.

It is worth mentioning that in sharp contrast with these two
problems, Connes's seminal work on injective factors \cite{connes1}
shows that if a group $G$ is ICC and amenable, $L(G)$ is isomorphic
to the unique hyperfinite $\II_1$ factor $R$.

\subsection{The group-measure space construction.} Let $G$ be a
countably infinite discrete group which acts in a measure preserving
way on a Borel probability space $(X,\mu)$. For each $g\in G$ and
$\zeta\in L^2(X,\mu)$ the formula
$$
\sigma_g(\zeta)(x)=\zeta(g^{-1}\cdot x)
$$
defines a unitary operator on $L^2(X,\mu)$.

We identify the Hilbert space $\mathcal H=L^2(G,L^2(X,\mu))$ with the Hilbert space of formal sums $\sum_{g\in G}\zeta_g\xi_g$, where the coefficients $\zeta_g$ are in
$L^2(X,\mu)$ and satisfy $\sum_g\|\zeta_g\|_{L^2(X,\mu)}^2<\infty$, and
$\xi_g$ are indeterminates indexed by the elements of $G$. The inner product on $\mathcal H$ is given by
$$
\langle \sum_{g\in G}\zeta_g(x)\xi_g, \sum_{g\in G}\zeta'_g(x)\xi_g\rangle=\sum_{g\in G} \langle\zeta_{g},\zeta'_{g}\rangle_{L^2(X,\mu)}.
$$
Both
$L^\infty(X,\mu)$ and $G$ act by left multiplication on $\mathcal H$
by the formulas
\begin{align*}
&f(\zeta_g(x)\xi_g)=((f(x)\zeta_g(x))\xi_g,\\
&u_h(\zeta_g(x)\xi_g)=\sigma_h(\zeta_g)(x)\xi_{hg},
\end{align*}
where $f\in L^\infty(X,\mu)$, $\zeta_g(x)\in L^2(X,\mu)$ and $g,h\in
G$. Thus if we denote by $\mathcal{FS}$ the set of finite sums,
$$
\mathcal{FS}=\{\sum_{g\in G}f_gu_g: \,f_g\in L^\infty(X,\mu),\,\, f_g=0, \text{
except for finitely many } g\},
$$
then each element in $\mathcal{FS}$ defines a bounded operator on
$\mathcal H$. Moreover, multiplication and involution in
$\mathcal{FS}$ satisfy the formulas
$$
(f_gu_g)(f_hu_h)=f_g\sigma_g(f_h)u_{gh}
$$
and
$$
(fu_g)^*=\sigma_{g^{-1}}(f^*)u_{g^{-1}}
$$
and so $\mathcal{FS}$ is a $*$-algebra. By definition, the {\it
group-measure space von Neumann algebra} is the weak operator
closure of $\mathcal{FS}$ on $\mathcal B(\mathcal H)$ and it is
denoted by $L^\infty(X,\mu)\rtimes_\sigma G$. The trace on
$\mathcal{FS}$, defined by
$$
\tau( \sum_{g\in G}f_gu_g)=\int_X f_e\,d\mu,
$$
extends to a faithful normal tracial state in
$L^\infty(X)\rtimes_{\sigma} G$ by the formula $\tau(T)=\langle
T(\xi_e),\xi_e\rangle$, where $e$ represents the identity of $G$.
Observe that $L^\infty(X,\mu)$ embeds into
$L^\infty(X)\rtimes_{\sigma} G$ via the map $f\mapsto fu_e$ and has
the property that its normalizer inside
${L^\infty(X)\rtimes_{\sigma} G}$,
\begin{align*}
&\mathcal N_{L^\infty(X)\rtimes_{\sigma} G}(L^\infty(X,\mu))=\\
&\{u\in \mathcal U(L^\infty(X)\rtimes_{\sigma} G):\,
uL^\infty(X,\mu)u^*=L^\infty(X,\mu)\}
\end{align*}
generates a weakly dense subalgebra in $L^\infty(X)\rtimes_{\sigma}
G$.

If $\sigma$ is a free action, then $L^\infty(X,\mu)$ is a MASA
(maximal abelian subalgebra) of $M$, in which case it is called a
{\it Cartan subalgebra} of $L^\infty(X)\rtimes_{\sigma} G$, (i.e. a
MASA with a weakly dense normalizer). If $\sigma$ is free then
$L^\infty(X)\rtimes_{\sigma} G$ is a factor (of type $\II_1$, since
$G$ is infinite) if and only if $\sigma$ is an ergodic action.

There is an important connection between the notion of {\it orbit
equivalence} and certain isomorphisms between group-measure space
von Neumann algebras. Recall that if $\sigma$ and $\tau$ are measure
preserving actions on standard Borel probability spaces $(X,\mu)$
and $(Y,\nu)$, respectively, of possibly different groups $G$ and
$H$, we say that $\sigma$ and $\tau$ are {\it orbit equivalent} if
there is a measure preserving bijection $\theta: X\to Y$ such that
$$
xE_\sigma x'\iff \theta(x)E_\tau \theta(x') \text{\ (a.e.)},
$$
i.e. if $\sigma$ and $\tau$ generate ``isomorphic'' orbit
equivalence relations $E_{\sigma}$ and $E_{\tau}$. Feldman and Moore
showed in \cite{fm} that two free ergodic measure preserving actions
$\sigma$ and $\tau$ are orbit equivalent if and only if their
corresponding inclusions of Cartan subalgebras $L^\infty(X)\subset
L^\infty(X)\rtimes_{\sigma} G$ and $L^\infty(Y)\subset
L^\infty(Y)\rtimes_{\tau} H$ are isomorphic. Thus the study of orbit
equivalence of measure preserving group actions can be translated
into a problem regarding inclusions of finite von Neumann algebras.

\section{Outline of the proofs}

\subsection{Theorem 1} Let $a:\SL(2,\Z)\actson\Z^2$ be the natural linear
action of $\SL(2,\Z)$ on $\Z^2$. Consider the natural measure
preserving ergodic a.e. free action $\sigma_0$ of $\SL(2,\Z)$ on
$X=\T^2$ equipped with the Haar measure $\mu$ and given
by
$$
\sigma_0(g)(\chi)(h)=\chi(a(g^{-1})(h)),
$$
where we identify $\T^2$ with the character group of $\Z^2$.
The matrices
$$
A=\left(\begin{array}{cc} 1 & 2\\ 0 & 1\end{array}\right),
B=\left(\begin{array}{cc} 1 & 0\\ 2 & 1\end{array}\right)
$$
generate a copy of $\F_2$ as a finite index subgroup of $\SL(2,\Z)$,
and so we get an action $\sigma:\F_2\actson\T^2$ by letting
$\sigma=\sigma_0|\F_2$. This action is still ergodic, see for instance \cite[\S
2]{torn1}. The group measure space factor
$L^\infty(X,\mu)\rtimes_{\sigma}\F_2$ was studied in detail by Sorin
Popa in \cite{popa}, where it was shown that $L^\infty(X,\mu)\subset
L^\infty(X,\mu)\rtimes_{\sigma}\F_2$ is a so-called $HT_s$ Cartan
subalgebra (see below), and that it is the unique $HT_s$ Cartan subalgebra in
$L^\infty(X,\mu)\rtimes_{\sigma}\F_2$ up to conjugation by a
unitary. In effect this means that the unitary conjugacy class of
$L^\infty(X,\mu)$ is definable inside of
$L^\infty(X,\mu)\rtimes_{\sigma}\F_2$ and depends only on the
isomorphism type of $L^\infty(X,\mu)\rtimes_{\sigma}\F_2$.

We refer to \cite[Definition 6.1]{popa} for the exact definition of
$HT_s$, but in short, the $H$ stands for Haagerup property, meaning
that $L^\infty(X,\mu)$ has the relative Haagerup property in
$L^\infty(X,\mu)\rtimes_{\sigma}\F_2$, and the $T$ means that it
also has the relative property (T), that is, $L^\infty(X,\mu)\subset
L^\infty(X,\mu)\rtimes_{\sigma}\F_2$ is a rigid inclusion in the
sense defined by Popa. These are von Neumann algebra generalizations
of the corresponding properties for discrete groups. Rather than
explaining the technical definition of the Haagerup property for
groups, we refer the reader to the monograph \cite{CCJJV}, and to
the recent survey paper \cite{pestov} in this journal for
applications and open questions regarding property H. Here we just
mention that amenable groups and the free groups $\F_n$ have the
Haagerup property \cite{HaagerupH}, and hence by \cite[Theorem
3.1]{popa} the inclusion $L^\infty(X,\mu)\subset
L^\infty(X,\mu)\rtimes_{\sigma}\F_2$ has the relative property H.
For groups, the property H and the property (T) are mutually
exclusive, in the sense that a discrete group that satisfies both
properties must be finite. A remarkable fact is that the inclusion
of groups $\Z^2\subset\Z^2\rtimes\SL(2,\Z)$ satisfies both the
relative property (T) and the property H. It is the combination of
deformation (i.e., the Haagerup property) and rigidity (i.e.,
property (T)), in particular the inclusion
$\Z^2\subset\Z^2\rtimes\SL(2,\Z)$, that is the engine behind Popa's
results in \cite{popa}, and in turn, the engine behind our results.

Returning to the proof, we now proceed as in \cite{torn1}: Let $a,b$
be generators for $\F_2$, and let $T_a,T_b\in\Aut(X,\mu)$ be the
measure preserving transformations corresponding to the action of
$a$ and $b$ according to $\sigma$. (Here $\Aut(X,\mu)$ denotes the
group of measure preserving transformations, equipped with Polish
group topology it inherits when naturally identified with a weakly
closed subgroup of the unitary group of $L^2(X,\mu)$.) It was shown
in \cite[\S 3]{torn1} that the set
\begin{multline*}
\Ext(\sigma)=\\
\{S\in\Aut(X,\mu): T_a,T_b \text{ and } S \text{ generate an a.e.
free action of } \F_3\}
\end{multline*}
forms a dense $G_\delta$ in $\Aut(X,\mu)$ and that if we, for $S\in
\Ext(\sigma)$, denote by $\sigma_S:\F_3\actson (X,\mu)$ the
resulting a.e. free ergodic $\F_3$-action, then the equivalence
relation
$$
S_1\sim_{oe} S_2\iff \sigma_{S_1} \text{ is orbit equivalent to }
\sigma_{S_2}
$$
has meagre classes and the set of transformations with dense $\sim_{oe}$-class is comeagre. It was pointed
out by Kechris in \cite[Theorem 17.1]{kechris} that this equivalence
relation is {\it generically $S_\infty$-ergodic}, meaning that if
$Y$ is a Polish $S_\infty$ space and $f:\Aut(X,\mu)\to Y$ is a Baire
measurable map which satisfies
$$
S_1\sim_{oe} S_2\implies (\exists g\in S_\infty) g\cdot
f(S_1)=f(S_2)
$$
then $f$ must be constant on a comeagre set. Since
$\simeq^{\Mod(\mathcal L)}$ is induced by a continuous $S_\infty$
action, this shows that $\sim_{oe}$ is not classifiable by countable
structures.

For $S\in\Ext(\sigma)$, let
$$
M_S=L^\infty(X,\mu)\rtimes_{\sigma_S}\F_3.
$$
The fact that $\F_3$ has the Haagerup property and that
$L^\infty(X,\mu)\rtimes_\sigma \F_2\subseteq
L^\infty(X,\mu)\rtimes_{\sigma_S}\F_3$ can  be seen to imply that
$L^\infty(X,\mu)$ is the unique (up to perturbation by a unitary)
$HT_s$ Cartan subalgebra of $L^\infty(X,\mu)\rtimes_{\sigma_S}\F_3$.

One now shows that the map $S\mapsto M_S$ is Borel. Further, if
$S\sim_{oe} S'$ then $M_S\simeq M_{S'}$ by Feldman and Moore's
Theorem. On the other hand, if $M_S\simeq M_{S'}$ then any
isomorphism $\varphi:M_S\to M_{S'}$ must, after possibly perturbing
it with a unitary, map $L^\infty(X,\mu)\subset M_S$ to
$L^\infty(X,\mu)\subset M_{S'}$. But then by Feldman-Moore, we must
have that $\sigma_S$ is orbit equivalent to $\sigma_{S'}$. Thus
$S\mapsto M_S$ provides a Borel reduction of $\sim_{oe}$ to
$\simeq^{\mathscr F_{\II_1}}$. Consequently, since $\sim_{oe}$ is
not classifiable by countable structures, neither is
$\simeq^{\mathscr F_{\II_1}}$.

\bigskip

The $\II_\infty$ and $\III_\lambda$ cases are consequences of the
$\II_1$ case, but this requires more sophisticated use of the
rigidity properties of the factors $M_S$ above. For the $\II_\infty$
case, one shows that the map
$$
S\mapsto M_S\otimes \mathcal B(l^2(\N)),
$$
where $\mathcal B(l^2(\N))$ denotes the bounded operators on
$l^2(\N)$, is a Borel reduction of $\sim_{oe}$ to
$\simeq^{\II_{\infty}}$. For the $\III_\lambda$ case, the map
$$
S\mapsto M_S\otimes R_\lambda
$$
provides a Borel reduction of $\sim_{oe}$ to
$\simeq^{\III_\lambda}$, where $R_\lambda$ is a (fixed) injective
factor of type $\III_\lambda$.

\subsection{Theorem 2.} Theorem 2 relies on another
deformation-rigidity result of Sorin Popa. Recall that if $G$ is a
countably infinite group then the (left) Bernoulli shift
$\beta:G\actson [0,1]^G$ is defined by
$$
\beta(g)(x)(h)=x(g^{-1}h).
$$
The Bernoulli shift is ergodic and preserves the product measure.

{\sc Theorem} (Popa, \cite[7.1]{popa2}). {\it Suppose $G_1$ and
$G_2$ are countably infinite discrete groups, $\beta_1$ and
$\beta_2$ are the corresponding Bernoulli shifts on
$X_1=[0,1]^{G_1}$ and $X_2=[0,1]^{G_2}$, respectively, and
$M_1=L^2(X_1)\rtimes_{\beta_1} G_1$ and
$M_2=L^2(X_2)\rtimes_{\beta_2} G_2$ are the corresponding
group-measure space $\II_1$ factors. Suppose further that $G_1$ and
$G_2$ are ICC groups having the relative property $(\pt)$ over an
infinite normal subgroup. Then $M_1\simeq M_2$ iff $G_1\simeq G_2$.}

The group $\SL(3,\Z)$ has property (T) outright (see \cite{behava})
and is ICC, and so any group of the form $H\times\SL(3,\Z)$, where
$H$ is ICC, satisfies the hypotheses of Popa's Theorem. Thus, to
prove Theorem 2, it suffices to show that if $\mathcal L$ is a
countable language, then $\simeq^{\Mod(\mathcal L)}$ is Borel
reducible to isomorphism of groups of the form $H\times\SL(3,\Z)$,
$H$ ICC, i.e. that isomorphism of groups of the form
$H\times\SL(3,\Z)$ is {\it Borel complete} for countable structures,
in the sense of \cite{friedstan}.

To this end, we modify a construction by Mekler, \cite{mekler}.
Mekler defines a notion of `nice graph', and proves (in effect) that
the isomorphism relation of countable connected nice graphs is Borel
complete for countable structures. Mekler then defines from a given
countable nice graph $\Gamma$ (and a prime $p$, which we shall keep
fixed here) a countable group $G(\Gamma)$, which we will call the
{\it Mekler group} of $\Gamma$, and shows that for nice graphs,
$\Gamma_1\simeq\Gamma_2$ iff $G(\Gamma_1)\simeq G(\Gamma_2)$. The
association $\Gamma\mapsto G(\Gamma)$ is Borel, and moreover, for
every graph automorphism of $\Gamma$ there is a corresponding group
automorphism of $G(\Gamma)$. However, the groups $G(\Gamma)$ are
generally {\it not} ICC.

To remedy this, we consider for each connected nice graph $\Gamma$
the nice graph $\Gamma_{\F_2}$, defined by
$$
(m,g)\Gamma_{\F_2} (n,h)\iff m\Gamma n\wedge g=h,
$$
consisting of $\F_2$ copies of $\Gamma$. ($\Gamma_{\F_2}$ is not
connected, but still nice.) Clearly, $\F_2$ acts by graph
automorphisms on $\Gamma_{\F_2}$. Going to the corresponding Mekler
group $G(\Gamma_{\F_2})$, we have a corresponding action of $\F_2$
by group automorphisms on $G(\Gamma_{\F_2})$. Thus we may form the
semi-direct product $G(\Gamma_{\F_2})\rtimes\F_2$. One now checks
that this groups is indeed ICC. Thus the group
$$
G_\Gamma=\SL(3,\Z)\times G(\Gamma_{\F_2})\rtimes\F_2,
$$
we obtain an ICC group with the relative property (T) over
$\SL(3,\Z)$. The argument is finished by arguing that $\SL(3,\Z)$,
as a subgroup of $G_\Gamma$, consists exactly of the elements of
$G_\Gamma$ which commutes with all elements of
$$
\{g\in G_\Gamma: (\exists\chi\in\Char(G_\Gamma))\ \chi(g)\neq 1\}.
$$
Supposing now that $G_{\Gamma^1}$ and $G_{\Gamma^2}$ are isomorphic,
it follows that $G(\Gamma^1_{\F_2})\rtimes\F_2$ is isomorphic to
$G(\Gamma^2_{\F_2})\rtimes\F_2$, from which it may in turn be
deduced that $G(\Gamma^1_{\F_2})$ is isomorphic to
$G(\Gamma^2_{\F_2})$. Then by Mekler's construction,
$\Gamma^1_{\F_2}\simeq \Gamma^2_{\F_2}$, so if $\Gamma^1$ and
$\Gamma^2$ are connected nice graphs then $\Gamma^1\simeq\Gamma^2$.
Thus the isomorphism relation of connected nice graphs is Borel
reducible to isomorphism of countable groups with the relative
property (T) over an infinite normal subgroup, which by Popa's
Theorem is all we needed to show.

\section{Some open problems} In this section we briefly discuss some
open problems related to the results stated above that we find may
be of interest to logicians.
\subsection{The Effros-Mar\'echal topology.} The space $\vN(\mathcal
H)$ has a natural Polish topology, called the {\it
Effros-Mar\'echal} topology. It is most easily defined as follows:
Let $L_1(\mathcal H)$ denote the unit ball in $\mathcal B(\mathcal
H)$, which is compact in the weak topology. Then the map
$$
M\mapsto M\cap L_1(\mathcal H)
$$
is 1-1, and so we may identify $\vN(\mathcal H)$ with a subset of
$K(L_1(\mathcal H))$, the space of compact subsets of $L_1(\mathcal
H)$. The Effros-Mar\'echal topology is the topology $\vN(\mathcal
H)$ inherits under this identification. It may be shown that
$\vN(\mathcal H)$ is Polish in this topology, see \cite[Theorem
2.8]{haagwin1}. The set of factors $\mathscr F$ forms a dense
$G_\delta$ set in $\vN(\mathcal H)$, see \cite[p. 402]{haagwin2}.
The most fundamental open problem seems to be:

\newtheorem{problem}{Problem}
\begin{problem}
Are the isomorphism classes in $\vN(\mathcal H)$ (equivalently,
$\mathscr F$) meagre? Are the unitary conjugacy classes meagre?
\end{problem}

\noindent If either part of problem 1 is answered in the
affirmative, the next natural question to ask is:

\begin{problem}
Does the unitary group act turbulently on $\vN(\mathcal H)$?
\end{problem}

Even though the subsets of $\II_1$, $\II_\infty$ or $\III$ factors
are not Polish in the Effros-Mar\'echal topology, they all form
Borel sets, and it is tempting to ask if one can find `natural'
topologies on these spaces in which Problem 1 and 2 would make
sense. We remark that by \cite[5.2.1]{beckerkechris}, it is possible
to find a Polish topology (with the same Borel structure) on the
subsets $\II_1$, $\II_\infty$ and $\III$ such that the conjugation
action of the unitary group becomes a continuous action, but by the
same token, \cite[5.1.6]{beckerkechris}, applying this too crudely
might make a conjugacy class clopen. Thus what we are really asking
is if these sets can be given Polish topologies where Problem 1 and
2 have affirmative answers.

It should be noted that Problem 1 is strongly related to the
so-called Connes embedding conjecture (see for instance \cite{pestov}) for separable von Neumann
algebras, which states that every separable type $\II_1$ factor can
be embedded into the ultrapower $R^\N/\mathcal U$, where $R$ is the
injective type $\II_1$ factor, and $\mathcal U$ is an ultrafilter
in $\N$. ($R^\N/\mathcal U$ is usually denoted $R^\omega$ in the von
Neumann algebra literature, since it is convention there to use $\omega$ to denote the ultrafilter.) Indeed, an affirmative answer to
Problem 1 is tantamount to refuting this conjecture, by the work of
Haagerup and Winsl\o w in \cite{haagwin2}. Namely, Haagerup and
Winsl\o w have shown that the Connes' embedding conjecture is
equivalent to the statement that the injective factors are dense in
$\mathscr F$. Since the set of injective factors is $G_\delta$,
Connes' embedding conjecture is equivalent to that the generic
element in $\mathcal F$ is injective. On the other hand, Haagerup
and Winsl\o w have also shown that the type $\III_1$ factors form a
dense $G_\delta$ subset of $\mathcal F$. Hence the Connes embedding
conjecture is equivalent to the assertion that the isomorphism class
of the (unique) injective type $\III_1$ factor forms a dense
$G_\delta$ set.

\subsection{$\itpf1$ factors and $T$-sets.} A factor $M$ is called an
ITPFI factor (short for Infinite Tensor Product of Factors of type I, also called an Araki-Woods factor), if it has the form
$$
M=\bigotimes_{k=1}^\infty (M_{n_k}(\C),\phi_k)
$$
where $M_{n_k}(\C)$ denotes the algebra of $n_k\times n_k$ matrices
and the $\phi_k$ are faithful normal states. (We refer the reader to
\cite[III.3.1]{blackadar} for the necessary basics regarding
infinite tensor products.). Among the ITPFI factors, the {\it Powers
factors} $R_\lambda$,  $0<\lambda<1$, are defined by taking $n_k=2$
for all $k$ and $\phi_k(x)=Tr(\rho_\lambda x)$ where
$$
\rho_\lambda=
\begin{pmatrix}
\frac{1}{1+\lambda}&0\\
0&\frac{\lambda}{1+\lambda}\\
\end{pmatrix}.
$$
Historically, the importance of the Powers factors is twofold: they
provided the first example of uncountably many non isomorphic von Neumann
factors (all of type $\III$) \cite{powers}. They were also the starting point of
the asymptotic analysis of factors carried on in the late sixties and early seventies that culminated
 with Connes classification of type $\III$ factors \cite{connes2}
and of injective factors \cite{connes1}. Since $\itpf1$ factors are
in particular injective factors, a corollary of Connes work is that
up to isomorphism there is only one $\itpf1$ factor of type
$\III_\lambda$, for each $\lambda\neq0$.  At the same time, Woods
proved \cite{woods} that the classification problem for ITPFI
factors is not smooth by showing that $E_0$ is Borel reducible to
isomorphism of ITPFI factors. (see \cite{woods2} for an historical
overview of $\itpf1$ factors and chapter \S 5 of \cite{connesNCG}
for an overview of Connes work.) Of course, the factors analyzed by
Woods in \cite{woods} are of type $\III_0$ and injective. In \S 4 of
our forthcoming paper \cite{asro} we show
\begin{theorem}
The isomorphism relation for injective factors of type $\III_0$ is
not classifiable by countable structures.
\end{theorem}

However the following remains open:
\begin{problem}
Are $\itpf1$ factors classifiable by countable structures?
\end{problem}

Woods uses an invariant $\rho(M)$, defined as
$$
\rho(M)=\{\lambda\in(0,1): M\otimes R_\lambda\simeq R_\lambda\}
$$
to distinguish the factors constructed there up to isomorphism.
(Here $R_\lambda$, $\lambda\in (0,1)$, denotes the
Powers factors, see \cite[III.3.1.7]{blackadar}). The invariant
$\rho$ has been replaced by the Connes invariant $T(M)$, called
the {\it $T$-set} of $M$, the general definition of which is rather
intricate. In the context of ITPFI factors, $T(M)$ is given by
$$
T(M)=\{t\in\R:\sum^{\infty}_{i=1}\big(1-|\sum_{k}
(\alpha_k^{(i)})^{1+it}|\big)<\infty\},
$$
where $\alpha_k^{(i)}$ denotes $k$th eigenvalue of $\phi_i$, see
\cite[III.4.6.9]{blackadar}. From this it can be deduced that the
$T$-set is a $K_\sigma$ subgroup of $\R$. It has been shown that all
countable subgroups of $\R$ and many uncountable subgroups are
realizable as $T$-sets of an $\itpf1$ factor (see \cite{gioskan}),
but the following seems to be open:

\begin{problem}
Is every $K_\sigma$ subgroup of $\R$ the $T$-set of some $\itpf1$
factor?
\end{problem}

The most natural approach to this problem would be to try to
construct from a given $K_\sigma$ subgroup $G\leq\R$ a corresponding
$\itpf1$ factor $M$ with $T(M)=G$. Can such a construction be
natural? More precisely, let
\begin{align*}
\mathscr S_{\sigma}(\R)=\{(K_n)\in K(\R)^\N: &(\forall n)
K_n=-K_n\wedge K_n\subseteq K_{n+1}\\
&\wedge K_n+K_n\subseteq K_{n+1} \}
\end{align*}
and let
$$
(K_n)\sim (K'_n)\iff \bigcup K_n=\bigcup K'_n.
$$
Then we can identify a $K_\sigma$ subgroup of $\R$ with an
equivalence class in $\mathscr S_{\sigma}(\R)/\sim$.

\newtheorem{problem5}[problem]{Problem}
\begin{problem5}
Is there a Borel $f:\mathscr S_{\sigma}(\R)\to\itpf1$ such that
$$
T(f(K_n))=\bigcup K_n
$$
and if $(K_n)\sim (K'_n)$ then $f(K_n)\simeq f(K'_n)$?
\end{problem5}

\subsection{Group von Neumann algebras vs. group-measure space von Neumann algebras.}

It is clear from the outline of the proof of Theorem 1 that what we
have really shown is that $\II_1$ factors that arise from the group
measure space construction are not classifiable by countable
structures. Even more specifically, we are dealing with those that
arise from an $\F_3$-action. The proof may be adapted to show that
for any $n\geq 2$, the group measure space von Neumann algebras
arising from an $\F_n$ action are not classifiable by countable
structures. However, in the light of the results of \cite{ikt}, it
is natural to ask:

\begin{problem}
If $G$ is a countably infinite non-amenable group, is it true that
the group measure space $\II_1$ factors arising from probability
measure preserving ergodic $G$-actions are not classifiable up to
isomorphism by countable structures?
\end{problem}

Our last problem is about the contrasting situation for group von
Neumann algebras. It is generally known that group von Neumann
algebras and group measure spaces von Neumann algebras can be rather
different. Indeed, one of the most striking applications of free probability theory
is Voiculescu's Theorem \cite{VoicCartan} stating that the group von Neumann algebras $L(\F_n)$ of free groups
on $n$ generators $n\geq 2$, don't have Cartan subalgebras, thus they are not group measure space von Neumann algebras. This result was generalized recently on \cite{PopaOzawa} without using free probability theory.

Recall the equivalence $\sim_{\vN}$ from our discussion of the group
von Neumann algebra (\S 2.1):
$$
G\sim_{\vN} H\iff L(G) \text{ is isomorphic to } L(H).
$$
In light of Theorem 1, it is natural to ask:

\begin{problem}
Is the isomorphism relation for group von Neumann algebras of
countable groups classifiable by countable structures? That is, is
$\sim_{\vN}$ classifiable by countable structures?
\end{problem}

One could consider Problem 7 more narrowly and ask if a ``weak''
version of Connes' conjecture, discussed in \S 2.1, is true: Is the
relation $\sim_{\vN}$ restricted to the class of ICC property (T)
groups classifiable by countable structures? A negative answer to
this would of course refute Connes' conjecture in a very strong way.

On the other hand one could, more broadly, ask if the classification
problem for group von Neumann algebras is as difficult as the one
for group-measure space von Neumann algebras. More precisely, what
is the relationship between $\sim_{\vN}$ and the isomorphism
relation for group-measure space factors in the Borel reducibility
hierarchy?

All of these questions are, to our knowledge, wide open and quite
interesting, since their solution may shed some light on Connes'
conjecture and the general relationship between a group and its
group von Neumann algebra.

\bibliographystyle{asl}
\bibliography{Bulletin_final}

\end{document}